\documentclass{article}
\usepackage{graphicx} 
\usepackage{tikz}
\usepackage{listings}
\usepackage{amsmath}
\usepackage{amssymb}
\title{Bridging Relations between SVD in Tensor Networks and Common Matrix Operations in Quantum Information Theory}
\author{John.TM.Campbell }
\date{Department of Engineering and Computer Science, Eolas, Maynnoth University , Republic of Ireland \\ \vspace{0.5cm} 3rd February 2024 }

\begin{document}

\maketitle

\section*{Abstract}
This document explores the connections between singular value decomposition (SVD) and various tensor operations in the context of quantum state calculations and tensor networks. The relationships between SVD and quantities such as trace, trace distance, Frobenius norm, and fidelity are examined, highlighting their similarities and differences. The intuition behind these mathematical connections is explained, providing insights into the underlying principles and applications of SVD and tensor operations in quantum physics and computational mathematics.

\section{Introduction}
Singular value decomposition (SVD) is a powerful matrix decomposition technique that has found numerous applications in various fields, including quantum physics and computational mathematics. In parallel, tensor operations play a crucial role in quantum state calculations and tensor network computations. This document aims to explore the connections between SVD and tensor operations, shedding light on their similarities and differences.

In tensor networks, the singular value decomposition (SVD) plays a crucial role in representing and manipulating quantum states. The SVD is a factorization of a matrix or tensor into three components: U, S, and V†, where U and V are unitary matrices and S is a diagonal matrix with non-negative entries called singular values.

The concept of Schmidt rank is specific to bipartite quantum systems and is related to the entanglement between the two subsystems. It refers to the number of non-zero singular values obtained from the SVD of the state's density matrix when the system is divided into two parts.

The quantum state of the bipartite system, lets call it $\psi$, can be said to be entangled if its Schmidt rank (i.e. number of singular values) is strictly greater than 1, and is not entangled otherwise. Hence Schmidt ranks for networks are an important product of SVD for networks to determine their "entanglement conductivity" as it were.

\section{Motivations}
Understanding the relationships between SVD and tensor operations is of great interest due to several motivations. Firstly, these connections provide a deeper insight into the mathematical foundations of quantum state calculations and tensor networks. By leveraging the principles of SVD, researchers can gain a better understanding of the underlying structures and properties of quantum states and tensors.

Secondly, the exploration of these connections offers practical benefits in computational mathematics. The use of SVD in computing quantities such as trace, trace distance, Frobenius norm, and fidelity provides efficient and accurate methods for evaluating these metrics in tensor networks. This line of reasoning can lead to improved algorithms and computational techniques for analyzing and manipulating tensors.

\section{Observed relations so far}

\begin{enumerate}

\item SVD and Trace:
\begin{itemize}
    \item The trace of a matrix can be obtained by summing the singular values (diagonal elements) of its SVD \cite{golub2012matrix}.
    \item For a tensor network, the trace of a tensor contraction can be computed by performing SVD on the resulting tensor and summing the singular values \cite{orus2014practical}.
\end{itemize}

\item SVD and Trace Distance:
\begin{itemize}
    \item The trace distance between two quantum states can be computed as the sum of the absolute differences between the singular values of their SVD \cite{nielsen2010quantum}.
    \item In tensor network calculations involving quantum states, the SVD can be used to compute the trace distance between two states by comparing their singular values \cite{orus2014practical}.
\end{itemize}

\item SVD and Frobenius Norm:
\begin{itemize}
    \item The Frobenius norm of a matrix is equal to the square root of the sum of the squares of its singular values \cite{golub2012matrix}.
    \item Similarly, in tensor networks, the SVD can be used to compute the Frobenius norm of a tensor by taking the square root of the sum of the squares of its singular values \cite{orus2014practical}.
\end{itemize}

\item SVD and Fidelity:
\begin{itemize}
    \item The fidelity between two quantum states can be computed as the square root of the sum of the squares of the singular values of their SVD \cite{nielsen2010quantum}.
    \item In tensor network calculations involving quantum states, the SVD can be used to compute the fidelity between two states by considering their singular values \cite{orus2014practical}.
\end{itemize}

\end{enumerate}

\section{intuition behind the maths}

1. SVD and Trace:
   For a matrix $A \in \mathbb{C}^{m \times n}$, its singular value decomposition (SVD) is given by $A = U \Sigma V^*$, where $U \in \mathbb{C}^{m \times m}$ and $V \in \mathbb{C}^{n \times n}$ are unitary matrices, and $\Sigma \in \mathbb{R}^{m \times n}$ is a diagonal matrix containing the singular values of $A$ on its diagonal.
   - The trace of a matrix $A$ can be computed as $\text{tr}(A) = \sum_{i=1}^{\min(m,n)} \sigma_i$, where $\sigma_i$ are the singular values of $A$.

2. SVD and Trace Distance:
   - The trace distance between two quantum states $\rho$ and $\sigma$ can be computed as $\text{Tr}|\rho - \sigma|$, where $|\cdot|$ denotes the absolute value and $\text{Tr}$ represents the trace operation.
   - In tensor network calculations involving quantum states, the SVD can be used to compute the trace distance between two states by comparing the singular values obtained from their respective SVDs.

3. SVD and Frobenius Norm:
   - The Frobenius norm of a matrix $A \in \mathbb{C}^{m \times n}$ is defined as $\|A\|_F = \sqrt{\sum_{i=1}^{m}\sum_{j=1}^{n} |a_{ij}|^2}$, where $a_{ij}$ represents the entries of $A$.
   - Similarly, in tensor networks, the SVD can be used to compute the Frobenius norm of a tensor by taking the square root of the sum of the squares of its singular values.

4. SVD and Fidelity:
   - The fidelity between two quantum states $\rho$ and $\sigma$ is defined as $F(\rho, \sigma) = \text{Tr}\sqrt{\sqrt{\rho}\sigma\sqrt{\rho}}$, where $\sqrt{\cdot}$ denotes the matrix square root.
   - In tensor network calculations involving quantum states, the SVD can be used to compute the fidelity between two states by considering the singular values obtained from their respective SVDs.

\begin{table}[h]
\centering
\caption{Differences and Similarities between SVD and Tensor Operations}
\resizebox{\textwidth}{!}{%
\begin{tabular}{|p{6cm}|p{6cm}|}
\hline
\textbf{Differences} & \textbf{Similarities} \\ \hline
SVD is a matrix decomposition technique & Both SVD and tensor operations are used in quantum state calculations and tensor networks \\ \hline
SVD decomposes a matrix into three separate matrices & Both SVD and tensor operations involve the use of singular values \\ \hline
Tensor operations involve various operations on tensors & SVD provides a mathematical framework for computing quantities in tensor networks \\ \hline
SVD is primarily used for analyzing and manipulating matrices & Both SVD and tensor operations can compute quantities like trace, trace distance, Frobenius norm, and fidelity \\ \hline
SVD involves unitary and diagonal matrices & Tensor operations can involve reshaping, tensor contractions, and other operations specific to tensors \\ \hline
\end{tabular}%
}
\end{table}

\section{Applications}
The understanding of the connections between SVD and tensor operations has several valuable applications. In quantum physics, these connections enable researchers to compute and compare quantities like trace distance and fidelity between quantum states. This information is crucial for analyzing the similarity or distinguishability of quantum systems and for measuring the accuracy of quantum state preparation and manipulation.

In the field of tensor networks, these connections allow for efficient computation of important metrics such as trace, Frobenius norm, and fidelity. These metrics play a vital role in assessing the quality and accuracy of tensor network approximations, enabling researchers to optimize tensor network calculations and improve the efficiency of computational simulations.

Furthermore, the insights gained from understanding the connections between SVD and tensor operations can be applied to various other domains involving matrix and tensor analysis. This includes applications in signal processing, image and video compression, machine learning, and data analysis, where SVD and tensor operations are commonly used for dimensionality reduction, feature extraction, and data compression.

The formula for SVD is as follows:
Singular Value Decomposition (SVD):
\begin{center}
    
\begin{align}
    A &= U \Sigma V^T \\
    &= \sum_{i=1}^{r} u_i \sigma_i v_i^T
    \label{eq:svd}
    \tag{[1]}
\end{align}

\vspace{0.5cm}

The Schmidt rank is closely related to the singular value decomposition (SVD). When we perform the SVD on a bipartite state, it essentially decomposes the state into a sum of rank-1 terms. These rank-1 terms are the product of one singular value (also known as Schmidt coefficient) and the corresponding column vectors from the left and right singular matrices.

In Equation [1], $A$ is decomposed into the matrices $U$, $\Sigma$, and $V^T$, where $r$ represents the rank of the matrix. Here, $U$ and $V$ are unitary matrices, and $\Sigma$ is a diagonal matrix containing the singular values (Schmidt coefficients) in non-increasing order. The rank of $A$ is equal to the number of non-zero singular values, which is also referred to as the Schmidt rank.

To obtain the Schmidt rank decomposition, we can rewrite the SVD equation in terms of the Schmidt coefficients:
\vspace{0.5cm}

Schmidt Rank Decomposition:
\begin{align}
    |\Psi\rangle &= \sum_{i=1}^{r} \sigma_i |u_i\rangle \otimes |v_i\rangle \\
    &= \sum_{i=1}^{r} \sqrt{\lambda_i} |i\rangle \otimes |i\rangle
    \label{eq:schmidt_rank_decomposition}
    \tag{[2]}
\end{align}
\end{center}

In Equation [2], $|\Psi\rangle$ represents a quantum state that is decomposed into the Schmidt coefficients $\sigma_i$, and the corresponding orthogonal states $|u_i\rangle$ and $|v_i\rangle$. The sum over $i$ runs from 1 to $r$, and $\lambda_i$ represents the eigenvalues associated with the Schmidt coefficients.

Here, $r$ is the Schmidt rank, and $u_i$ and $v_i$ are the $i$th columns of $U$ and $V$ respectively. Each term $\sigma_i u_i v_i^T$ represents a rank-1 contribution to the overall state. Therefore, the Schmidt rank decomposition expresses the bipartite state $|\Psi\rangle$ as a sum of rank-1 terms, where the Schmidt coefficients $\sigma_i$ and the corresponding vectors $u_i$ and $v_i$ capture the entanglement structure of the state.

\vspace{1cm}

Some code listing in python is shown below for common tensor operations. In this code, some tensor operations are shown along with the SVD operation `np.linalg.svd(A)`, performed on the matrix `A`, and the resulting matrices `U`, `S`, and `VT` represent the left singular vectors, singular values, and right singular vectors, respectively:
\begin{lstlisting}[language=Python, caption={Tensor operations code}]
import numpy as np

# Perform tensor dot product
result = np.tensordot(tensor1, tensor2, axes=1)

# Perform tensor contraction
result = np.einsum('ijk,jkl->il', tensor1, tensor2)

# Perform SVD decomposition
U, S, VT = np.linalg.svd(A)

\end{lstlisting}

For furthering applications we also seek to borrow some concepts that bear similarities to the products of SVD, i.e. the Schmidt rank, from other fields to provide benchmarks of quantum networks. In this we are inspired to treat the quantum networks as if they are many-body quantum systems seen in the field of particle physics and quantum chemistry.

For instance, in quantum chemistry there is a very useful concept known as the Slater rank.

The Slater rank can be a measure of entanglement or correlation in a many-body wave function. It is defined as the minimum number of determinants required to represent the wave function accurately. A determinant in this context refers to a Slater determinant, which is a specific type of wave function used to describe the behavior of fermions.

The Slater rank can be related to the Schmidt rank in bipartite systems. In the case of a bipartite system, the Schmidt rank corresponds to the number of singular values obtained from the SVD of the density matrix when the system is divided into two parts. Similarly, the Slater rank corresponds to the number of determinants needed to accurately represent the many-body wave function.

In tensor network methods, such recent developments in the Multi-Scale Entanglement Renormalization Ansatz (MERA) \cite{vidal2008quantum},  method has given new
impetus to its application for strongly correlated systems. In the developments of MERA The Slater rank plays a crucial role in approximating the wave function using tensor networks \cite{entopabinitio2015}. By limiting the Slater rank, one can efficiently represent highly entangled states and capture the relevant physics of the system.

\section{Tensor network diagrams}

adopting the use of Tensor network diagrams can be a useful way to saliently understand some of the operations discussed here and develop a bit more intuition about their use in quantum networks.

This is a maturing branch of mathematical representation and a lot can be gleaned from literature while at the same time a particular flavour can be added by our own exploration into the utility of tensor and matrix operations that are being developed into practical benchmarks for quantum networks.

We employ the TIKz graph package to start to draw some of the established tensor network diagrams and attempt to draw a new diagram for representing the ranks discussed.

\begin{figure}[hbt!]
    \centering
    \includegraphics[width=1\textwidth]{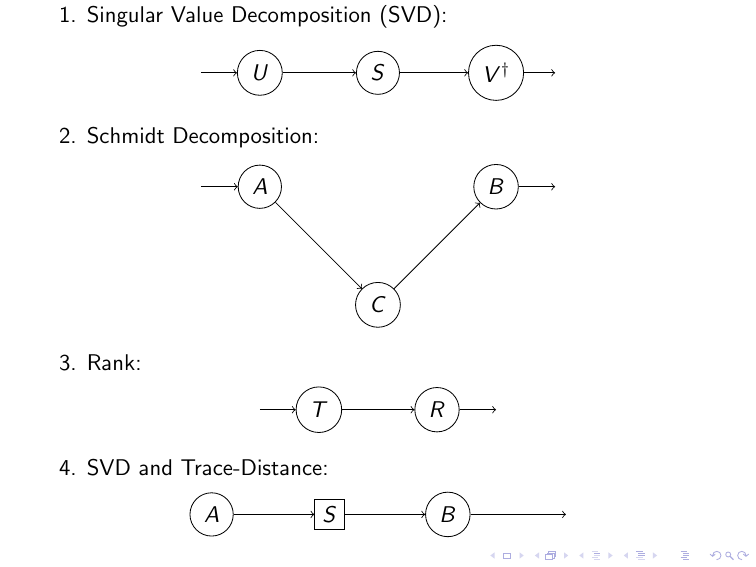}
\end{figure}

\vspace{10cm}

\section{Conclusions}

In conclusion, the exploration of the connections between SVD and tensor operations provides valuable insights into the mathematical foundations and practical applications of quantum state calculations and tensor networks. This line of reasoning has the potential to enhance our understanding of quantum physics, improve computational techniques, and find applications in a wide range of domains involving matrix and tensor analysis.

\end{document}